\documentclass[12pt]{amsart}
\usepackage{fancyhdr}
\usepackage{stmaryrd}
\usepackage{calrsfs}

\usepackage{amsmath}
\usepackage[nospace,noadjust]{cite}
\usepackage{amsfonts}
\usepackage{amssymb,enumerate}
\usepackage{amsthm}
\usepackage{cite}
\usepackage{comment}
\usepackage{color}
\usepackage[all]{xy}
\usepackage{hyperref}
\usepackage{lineno}
\usepackage{stmaryrd}

\usepackage{mathtools}
\usepackage{bm}
\usepackage{graphicx}
\usepackage{psfrag}
\usepackage{url}

\usepackage{fancyhdr}
\usepackage{tikz-cd}

\usepackage{dirtytalk}
\usepackage{float}
\usepackage{mathrsfs}
\usepackage{enumitem}

\newtheorem*{maintheorem*}{Main Theorem}
\newtheorem{theorem}{Theorem}[section]
\newtheorem{prop}[theorem]{Proposition}

\newtheorem{conj}[theorem]{Conjecture}
\newtheorem{lemma}[theorem]{Lemma}
\newtheorem{cor}[theorem]{Corollary}

\theoremstyle{definition}
\newtheorem{definition}[theorem]{Definition}
\newtheorem{remark}[theorem]{Remark}
\newtheorem{example}[theorem]{Example}
\numberwithin{equation}{section}

\newcommand{\nn}{\mathbb{N}}

\newcommand{\zz}{\mathbb{Z}}

\newcommand{\A}{\mathcal{A}}
\newcommand{\Z}{\mathsf{Z}}
\newcommand{\dd}{\mathsf{d}}
\newcommand{\cc}{\mathsf{c}}

\newcommand{\red}{\text{red}}
\newcommand{\bul}{\operatorname{bul}}

\hypersetup{
	pdftitle={FD},
	colorlinks=true,
	linkcolor=blue,
	citecolor=cyan,
	urlcolor=wine
}


\begin{document}

	\mbox{}
	\title{On the factorization invariants of arithmetical congruence monoids}
	\author{Scott T. Chapman}
	\author{Caroline Liu}
	\author{Annabel Ma}
	\author{Andrew Zhang}

\maketitle

\begin{abstract}
In this paper, we study various factorization invariants of arithmetical congruence monoids. The invariants we investigate are the catenary degree, a measure of the maximum distance between any two factorizations of the same element, the length density, which describes the distribution of the factorization lengths of an element, and the omega primality, which measures how far an element is from being prime.
\end{abstract}

\bigskip
\section{Introduction}
\label{sec:intro}

The fundamental theorem of arithmetic states that, for each integer $n$ greater than $1$, there is a unique factorization of $n$ into primes (up to permutation and multiplication by units). Yet, it is well known that this property does not hold for other algebraic structures such as rings of algebraic integers. This phenomenon of non-unique factorization led to Dedekind's ideal theory and Kroneker's divisor theory in the development of algebraic number theory during the 19th century.

Arithmetical congruence monoids (ACMs) are arithmetic progressions that are closed under multiplication. Specifically, we have that an ACM is a monoid of the form
\[
	M_{a,b} = \{a, a + b, a+2b, a+3b, \dots\} \cup \{1\} = (a+b\mathbb{N}_0) \cup \{1\}
\] 
for $a, b \in \mathbb{Z}$ such that $0 < a \leq b$ and $a^2 \equiv a \pmod b$. David Hilbert famously used these monoids as a pedagogical tool to demonstrate the necessity of proving the unique factorization property of the integers to his students (see ~\cite{davenport1999higher}). ACMs can exhibit both unique and nonunique factorization of elements. Consider the following examples. 

\begin{example}\label{ex: hilbert}
We consider the Hilbert monoid, defined as $$\textbf{H} = \{1+4k \ | \ k \in \mathbb{N}_0\} = 1+4\mathbb{N}_0.$$
Note that $693 = 21 \cdot 33 = 9 \cdot 77$. Additionally, we have $21 = 3 \cdot 7$, $33 = 3 \cdot 11$, $9 = 3 \cdot 3$, and $77 = 7 \cdot 11$. These factorizations in $\mathbb{N}$ imply that $21,$ $33,$ $9,$ and $77$ are all irreducibles in $\textbf{H}$. Thus, the monoid $\textbf{H}$ provides an example of an ACM displaying non-unique factorization.
\end{example}

\begin{example}\label{ex: natural numbers}
The set of natural numbers $\mathbb{N}$ is an ACM, namely $M_{1,1}.$ In contrast with the Hilbert monoid, $\mathbb{N}$ displays unique factorization into primes under the Fundamental Theorem of Arithmetic. We call such a monoid a \textit{unique factorization monoid (UFM)}.
\end{example}

The structure of ACMs is surprisingly complex. In fact, there exist various open problems regarding finer measures of the factorization invariants of ACMs (see, for example, \cite[Open Question~4.6 and Open Question~4.18]{BC14}). The purpose of the present article is to study the omega primality, length density, and catenary degree in the context of ACMs, which are values quantifying how far ACMs are from being UFMs. More specifically, the omega primality function measures how far elements of an ACM are from being prime, the length density measures the sparseness of factorization lengths within an ACM, and the catenary degree uses a notion of distance to bound how different factorizations of the same element can get. While other factorization properties of ACMs have been considered before~\cite{banister2007arithmetic, BCS08, banister2009theorem, hartzer2016periodicity, chapman2010elasticity}, the omega primality, length density, and catenary degree have not been researched as much in this context. 

Our paper is structured as follows. In Section~\ref{sec:background}, we review some of the standard notation and terminology we shall be using throughout the paper. In Section~\ref{sec: omega primality}, we compute a closed form of the omega primality of all ACMs, which depends on the powers of the primes dividing $a$ and $b$ in $M_{a,b}.$ In Section~\ref{sec: length density} we provide a closed form of the length density of a regular ACM $M_{1,b}$ based on $\phi(b)$. We also compute the length density of local singular ACMs. A conjecture regarding the closed form of the length density of global singular ACMs is also made. Finally, in Section~\ref{sec: catenary degree}, we compute the catenary degree of local singular ACMs by splitting this class into three cases. For a local singular ACM $M_{a,b}$ with $\gcd(a,b) = p^\alpha,$ we provide an explicit formula depending on $\alpha$ and the least power of $p$ within the monoid. We also provide conjectured closed forms for the catenary degree of global singular ACMs. 

\bigskip
\section{Preliminaries}
\label{sec:background}

Throughout this paper, a \emph{monoid} is defined to be a semigroup with identity that is cancellative and commutative. Unless otherwise specified, we will use multiplicative notation for monoids. Let $M$ be a monoid with identity~$1$. We set $M^{\bullet} \coloneqq M \setminus \{1\}$, and we let $M^{\times}$ denote the group of units (i.e., invertible elements) of~$M$. In addition, we let $M_{\red}$ denote the quotient $M/ M^{\times}$, which is also a monoid. The monoid $M$ is \emph{reduced} provided that $M^{\times}$ is the trivial group, in which case we naturally identify $M_{\red}$ with $M$. 

For $b,c \in M$, we say that $b$ \emph{divides} $c$ \emph{in} $M$ if there exists $b' \in M$ such that $c = b b'$, in which case we write $b \mid_M c$. Two elements $b,c \in M$ are \emph{associates} if $b \mid_M c$ and $c \mid_M b$. An element $a \in M \setminus M^{\times}$ is an \emph{atom} if for all $b,c \in M$ the equality $a = bc$ implies that either $b \in M^{\times}$ or $c \in M^{\times}$. On the other hand, an element $a \in M \setminus M^{\times}$ is \textit{prime} if $a \mid_M bc$ implies that either $a \mid_M b$ or $a \mid_M c$. We let $\mathcal{A}(M)$ denote the set of all atoms of $M$. The monoid $M$ is \emph{atomic} if each element in $M \setminus M^{\times}$ can be written as a (finite) product of atoms. One can readily check that~$M$ is atomic if and only if $M_{\red}$ is atomic. 

Assume now that $M$ is atomic. We let $\mathsf{Z}(M)$ denote the free (commutative) monoid on $\mathcal{A}(M_{\red})$. The elements of $\mathsf{Z}(M)$ are \emph{factorizations}, and if $z = a_1 \cdots a_\ell \in \mathsf{Z}(M)$ for $a_1, \ldots, a_\ell \in \mathcal{A}(M_{\red})$, then $\ell$ is the \emph{length} of $z$, which is denoted by $|z|$. Let $\pi \colon \mathsf{Z}(M) \to M_{\red}$ be the unique monoid homomorphism satisfying that $\pi(a) = a$ for all $a \in\mathcal{A}(M_{\red})$. For each $b \in M$, the following sets associated to $b$ are fundamental in the study of factorization theory:
\begin{equation} \label{eq:sets of factorizations/lengths}
	\mathsf{Z}_M(b) \coloneqq \pi^{-1}(b \mathcal{U}(M)) \subseteq \mathsf{Z}(M) \hspace{0.6 cm}\text{ and } \hspace{0.6 cm}\mathsf{L}_M(b) \coloneqq \{|z| : z \in\mathsf{Z}_M(b)\} \subseteq \nn_0.
\end{equation}
We drop the subscript $M$ in~\eqref{eq:sets of factorizations/lengths} whenever the monoid is clear from the context. 

In \cite{BC14}, the authors show that all ACMs fall into one of three mutually exclusive classes: regular, local singular, and global singular. A \textit{regular} ACM is an ACM of the form $M_{1,b}.$ These ACMs exist for all positive integers $b,$ as $a^2 = a$ implies that $a^2 \equiv a \pmod b$ for $a = 1$.

\begin{definition}
    Let $M$ be a monoid. A \textit{divisor theory} for $M$ is a free commutative monoid $\mathcal{F}(P)$ and a monoid homomorphism $\sigma : M \rightarrow \mathcal{F}(P)$ satisfying the following properties.
\begin{enumerate}
    \item $\sigma(u) = 1$ for any $u \in M^{\times}$;
    
    \item $\sigma(u) \neq 1$ for any $u \not\in M^{\times}$; 
    
    \item for any nonunits $x, y \in M,$ $\sigma(x) \mid_M \sigma(y)$ implies $x \mid_M y$; 
    
    \item for every $p \in P,$ there is a finite subset $X \subseteq M$ such that $p = \gcd(\sigma(X))$. 
\end{enumerate}
A monoid $M$ with a divisor theory is called a \textit{Krull monoid}. 
\end{definition}

We can use this idea to classify regular ACMs, in a theorem first shown in \cite{H-K91}.

\begin{theorem}\label{thm:krull}
Let $P = \{p \in \mathbb{N} \mid p \text{ is prime and } \gcd(p,b) =1\}.$ The free monoid $\mathcal{F}(P) \leq (\mathbb{N}, \times)$ and the homomorphism $\iota : M_{1, b} \rightarrow \mathcal{F}(P)$ form a divisor theory for $M_{1,b}.$ Thus $M_{1,b}$ is Krull. 
\end{theorem}

\textit{Singular} ACMs are monoids $M_{a,b}$ such that $\gcd(a,b) \neq 1$. The factorization structure of a singular ACM depends on $\gcd(a,b)$, which we will call $d$. We then set $f = b/d.$ Note that $d = 1$ if and only if $a = 1$. Then, an ACM $M_{a,b}$ is singular if and only if $a \neq 1$. Singular ACMs are divided into two classes based on $d$: \textit{local} ACMs have $d$ a power of a \textit{prime}, and global ACMs have $d$ divisible by more than one prime. 

The following two theorems introduced in~\cite{BCS08} create a criterion for inclusion of elements in a singular ACM and in the set of irreducibles of a singular ACM, which we will use extensively in later sections when computing our respective factorization invariants. 
\begin{theorem}\label{Theorem: p^beta 1 mod b'}
For a singular ACM $M_{a,b},$ we have $x \in M_{a,b}$ if and only if $x \equiv 1 \pmod b$.
\end{theorem}

\begin{theorem}\label{thm: division inclusion}
	Let $x, y \in M_{a,b}$ be such that $x, y \neq 1$ and $y \mid_{\mathbb{N}} x.$ Set $d = \gcd(a,b).$
	\begin{enumerate}
		\item If $d \mid_{\mathbb{N}} (x/y)$ then $x/y \in M_{a,b}$
		\item If $x \in \A(M_{a,b})$ then $y \in \A(M_{a,b}).$
	\end{enumerate}
\end{theorem}

Now, let $M$ be a monoid and define the \textit{delta set} of $x \in M$ as follows. 
\begin{definition}
Let $\mathsf{L}(x) = \{n_1, n_2, \dots n_k\}$ where $n_1 < n_2 < \dots < n_k$. Then $\Delta(x) = \{n_{i+1} - n_{i}| 1 \leq i < k \}$. Furthermore, 
\[
	\displaystyle{\Delta(M) = \bigcup_{x \in M} \Delta(x)}.
\]
\end{definition}
Finally, for a local singular ACM $M$, we know that $d = p^\alpha$ where $p$ is prime. We then set $\beta$ to be the smallest number such that $p^{\beta} \in M.$ 

Now that we have established some basic notation regarding ACMs, we will introduce the factorization invariants that we will discuss in this paper. 

The \emph{omega primality function} measures how far a nonunit element of a monoid is from being prime. First introduced in~\cite{geroldinger1997chains}, it is defined as follows.

\begin{definition}
\label{opf def}
    For some $x \in M$, we have $\omega_M(x) = m$, or simply $\omega(x) = m$ if the monoid in question is clear; if $m$ is the smallest positive integer such that if $x \mid \prod_{i = 1}^k a_i$ for $a_i \in M$ and $k > m$, there exists some proper subset $S \subset \{1, 2, \dots , k\}$ such that $x \mid \prod_{i \in S}a_i$. If no $m$ satisfying this condition exists, we let $\omega(x) = \infty$. 
\end{definition}

\begin{example}[{\cite[Example~2.5]{OP2017}}]
    In the monoid $(\nn, \times)$, we have that $\omega(p) = 1$ for any prime number $p$, and if $n = p_1p_2\cdots p_k$ for primes $p_1, p_2, \dots, p_k \in \nn$, then we have that $\omega(n) = k$. This is because if $n \mid a_1a_2\cdots a_\ell$ for $a_1, a_2, \dots, a_\ell \in \nn$, we must have that the product $p_1p_2\cdots p_k$ appears somewhere within the prime factorization of $a_1a_2\cdots a_\ell$. In the worst case scenario, each prime within the prime factorization of $n$ appears in a different $a_i$, so $\omega(n) \leq k$. And because $n = p_1p_2\cdots p_k$, we find that $\omega(n) = k$.
\end{example}

Definition~\ref{opf def} then prompts the definition of bullets. 
\begin{definition}
    For a monoid $M$ and $x \in M$, a \textit{bullet} of $x$ is a product $a_1a_2\cdots a_k$ where $a_i \in \mathcal{A}(M)$ for $i \in \{1, 2, \dots , k\}$ such that $x \mid a_1a_2\cdots a_k$ and $x$ does not divide the product of any proper subset of $\{a_1, a_2, \dots , a_k\}$. We denote $\bul(x)$ to be the set of bullets of $x$. 
\end{definition}

Proposition 2.10 of \cite{geroldinger1997chains} then defines the $\omega$-primality function in terms of bullets and provides a proof to show that it is analogous to Definition~\ref{opf def}. 
\begin{prop}
    For each $x \in M$ for a monoid $M$, we have that $$\omega(x) = \sup\{k\mid a_1a_2\cdots a_k \in \bul(x), a_i \in \mathcal(M)\}.$$
\end{prop}

Using both definitions of bullets, \cite[Proposition~2.12]{geroldinger1997chains} provides some basic properties of the $\omega$-primality function.
\begin{prop}
    In a commutative, cancellative, and atomic monoid $M$, the following statements hold.
    \begin{enumerate}
       \item The set $\{\omega(x)\mid x\in M\}$ is unbounded.
        \item For all $x, y \in M$, $\omega(xy) \leq \omega(x) + \omega(y)$. We then call $\omega$ \text{subadditive}.
        \item For some prime $p \in M$, we have $\omega(xp) = \omega(x) + 1$ for all $x \in M$.
   \end{enumerate}
\end{prop}

Length density, a concept introduced in ~\cite{chapman2022length}, is a novel tool used to analyze the sparseness of the distribution of the lengths of factorizations in monoids. We denote $\ell(x) = \min\mathsf{L}(x)$ and $L(x) = \max \mathsf{L}(x)$. We set L$^{\Delta}(x) \coloneqq L(x) - \ell(x).$ The set of elements $x \in M$ such that L$^{\Delta}(x) \neq 0$ is the \textit{length ideal} of $M$, denoted by $M^{LI}$. 
\begin{definition}
For every $x \in M^{LI}$ we define the \textit{length density of $x$} as $$\text{LD} (x) = \frac{|\mathsf{L}(x)|-1}{\textup{L}^{\Delta}(x)}.$$ Furthermore, the \textit{length density of $M$} is defined as 
$$
\text{LD}(M) = \text{inf\{LD}(x) \text{ } | \text{ } x \in M^{LI} \}.     
$$
\end{definition}

The catenary degree uses the idea of distance between factorizations of the same element to measure how close to a unique factorization monoid an ACM is. The \emph{distance function} is a metric that parameterizes ``how close" irreducible factorizations of an element of a monoid are. 

\begin{definition}
   For an atomic monoid $M$ and $x\in M,$ two factorizations $z_1, z_2 \in \Z(x)$ can be written as 
   \begin{align*}
       z_1 = \alpha_1\dots\alpha_j\gamma_1\dots\gamma_{\ell} \text{ and }z_2 = \beta_1\dots\beta_k\gamma_1\dots\gamma_{\ell},
   \end{align*}
   where $j, k, \ell \in \mathbb{N}_0$ and $\alpha_i, \beta_i, \gamma_i \in \A(M)$ such that $\{\alpha_1\dots\alpha_j\} \cap \{\beta_1\dots\beta_k\} = \emptyset$. Then, $\dd(z_1, z_2) = \max\{j, k\} \in \mathbb{N}_0$ is called the \textit{distance} between $z_1$ and $z_2$. 
   \end{definition}

We use the distance function to define a \textit{chain of factorizations}.

\begin{definition}
    Let $M$ be an atomic monoid and $x \in M.$ A sequence of factorizations $z_0, z_1, \dots, z_t$, where each $z_i \in \mathsf{Z}(x)$ is called a \textit{chain of factorizations} of $x.$ For each $1 \leq i \leq t,$ the $i$th link of the chain are the factorizations $z_{i-1},z_i.$ The \textit{length} of the $i$th link is $\mathfrak{d}_i = \dd(z_{i-1}, z_{i})$. 

\end{definition}

\begin{definition}
    Let $M$ be an atomic monoid with $x \in M.$ Let $N$ be a positive integer. A chain of factorizations $z_0, z_1, \dots, z_t$ in $\Z(x)$ is called an \textit{$N$-chain} if each distance $\mathfrak{d}_i \leq N$ for $i \in \{1,\dots, t\}.$
\end{definition}

We can now define the catenary degree.

\begin{definition}
     Let $M$ be an atomic monoid, and let $x \in M.$ The \textit{catenary degree of $x$} is defined as 
    $$\cc(x) = \text{min}\{N \mid \text{ there exists an $N$-chain between any } z_1, z_2 \in \Z(x)\}.$$
    We define the \textit{catenary degree of $M$} to be 
    $$\cc(M) = \sup\{\cc(x) \mid x \in M\}.$$
\end{definition}

\bigskip
\section{Omega Primality}
\label{sec: omega primality}

In this section, we provide closed formulas to compute the omega primality of the elements of an ACM. To do so, we first consider regular ACMs and then singular ACMs.

\begin{prop} \label{prop: omega primality regular case}
	Let $M_{1,b}$ be a regular ACM for some $b \in \nn_{>1}$, and let $x$ be an element of $M_{1,b}$. Then $\omega(x) = \sum_{i = 1}^{n} e_i$, where $\prod_{i = 1}^{n} p_i^{e_i}$ is the factorization of $x$ in $(\zz, \cdot)$.
\end{prop}

\begin{proof}
	First, let us show that $\omega(x) \leq \sum_{i = 1}^{n} e_i$. Observe that, for elements $c,d \in M_{1,b}$, we have that $c \mid_{M_{1,b}} d$ if and only if $c \mid_{\zz} d$. Indeed, if there exists $k \in \zz$ such that $ck = d$ then it is easy to see that $k \equiv 1 \pmod b$. Now let $z = a_1 \cdots a_m$ be an arbitrary bullet of $x$ and define $\llbracket a, b \rrbracket$ to be the positive integers between $a, b$, inclusive. For each $j \in \llbracket 1,m \rrbracket$, there exists $i \in \llbracket 1,n \rrbracket$ such that $\mathsf{v}_{p_i}(a_1\cdots \overline{a_j}\cdots a_m) < \mathsf{v}_{p_i}(x) = e_i$. Consequently, we have that $m \leq \sum_{i = 1}^{n} e_i$ which, in turn, implies that the inequality $\omega(x) \leq \sum_{i = 1}^{n} e_i$ holds by virtue of \cite[Proposition~2.10]{OP2017}.
	
	Fix $i \in \llbracket 1,n \rrbracket$. Note that $p_i$ and $b$ are relatively prime positive integers. Dirichlet's Theorem states that if $\gcd(a,d) = 1,$ there exists infinitely many primes of the form $a+nd.$ So, there are infinitely many primes of the form $p_i + kb.$ Hence there exist distinct $k_1, \ldots, k_{e_i} \in \nn$ such that $q_j = p_i + k_jb$ is a prime number for every $j \in \llbracket 1, e_i \rrbracket$. Without loss of generality, assume that $q_j \nmid_{\zz} x$. Set $x_j^i \coloneqq p_iq_j^{\varphi(b) - 1}$ for $j \in \llbracket 1, e_i \rrbracket$, where $\varphi$ denotes the Euler's totient function; these are all elements of $M_{1,b}$ by Euler's Theorem stating that for any modulus $n$ and any integer $a$ coprime to $n$, one has $a^{\varphi(n)} \equiv 1 \pmod n$. Now consider the product
	\[
		z = \prod_{\substack{i \in \llbracket 1,n \rrbracket \\ j \in \llbracket 1,e_i \rrbracket }} \!\!x_j^i.
	\]
	There is no loss in assuming that $\gcd(x_j^i, x_{j'}^{i'})$ is either $1$ or a prime number provided that either $i \neq i'$ or $j \neq j'$. We may also assume that $x_j^i \in \mathcal{A}(M_{1,b})$, which implies that $z$ is a factorization of $x$. It is not hard to see that $x \mid_{M_{1,b}} z$. Moreover, we have that $x \nmid_{M_{1,b}} z(x_j^i)^{-1}$ for any $i \in \llbracket 1,n \rrbracket$ and any $j \in \llbracket 1,e_i \rrbracket$. Consequently, $z$ is a bullet of $x$ with $|z| = \sum_{i = 1}^{n} e_i$, which concludes our proof.
\end{proof}
We now consider singular ACMs.

\begin{prop} \label{prop: omega primality singular case}
Let $M_{a,b}$ be a singular ACM for some $a, b \in \nn_{>1}$ such that $a \leq b$, and let $\gcd(a,b) = \prod_{i = 1}^m q_i^{r_i}$ with $q_i \in \mathbb{P}$ and $r_i \in \nn$ for every $i \in \llbracket 1,m \rrbracket$. Let $x \in M_{a,b}$ such that
	\[
		x = \left(\prod_{i = 1}^{m} q_i^{r_i + e_i}\right) \left(\prod_{j = 1}^{n} p_j^{s_j}\right)
	\]
	with $p_j \in \mathbb{P}$ for all $j \in \llbracket 1,n \rrbracket$ and $p_j \neq q_i$ for any $j \in \llbracket 1,n \rrbracket$ and any $i \in \llbracket 1,m \rrbracket$. Thus,
	\[
		\omega(x) = \max\left\{ 1 + \left\lceil \frac{r_1 + e_1}{r_1} \right\rceil, \ldots, 1 + \left\lceil \frac{r_m + e_m}{r_m} \right\rceil, \sum_{j = 1}^{n} s_j \right\}.
	\]
\end{prop}

\begin{proof}
	Set $d \coloneqq \gcd(a,b) = \prod_{i = 1}^m q_i^{r_i}$. Observe that, for elements $c,e \in M_{a,b}$, we have that $c \mid_{M_{a,b}} e$ if and only if $c \mid_{\zz} e$. Indeed, if there exists $k \in \zz_{\geq 2}$ such that $ck = e$, then for $n \in \nn$ we have
	\[
		k \equiv \frac{c}{e} \equiv nd \equiv nd \cdot a \equiv \frac{c}{e} \cdot e \equiv a \pmod b. 
	\]
	Let $z = a_1 \cdots a_m$ be an arbitrary bullet of $x$. Reasoning as in the first paragraph of the proof of Proposition~\ref{prop: omega primality regular case}, it is easy to show that $\omega(x) < \infty$. Consequently, we may assume that $z$ is a bullet of maximal length. Next we show that $xd \mid \pi(z)$. Suppose that $x = a_1 \cdots a_m$. Now set $a_1' \coloneqq a_1q$ and $a_{m + 1}' \coloneqq aq$ for some prime $q$ satisfying that $q \nmid x$ and $q \equiv 1 \pmod b$. For $i \in \llbracket 2,m \rrbracket$, set $a_i' \coloneqq a_i$. It is easy to see that $z' = a_1' \cdots a_{m + 1}'$ is a bullet of $x$ with length bigger than $m$, a contradiction. Hence $xd \mid \pi(z)$. Since $d \mid_{M_{a,b}} a_i$ for each $i \in \llbracket 1,m \rrbracket$ and $xd \mid \pi(z)$, the inequality 
	\begin{equation} \label{eq: upper bound}
		m \leq \max\left\{ 1 + \left\lceil \frac{r_1 + e_1}{r_1} \right\rceil, \ldots, 1 + \left\lceil \frac{r_m + e_m}{r_m} \right\rceil, \sum_{j = 1}^{n} s_j \right\}
	\end{equation}
	holds.
	
	Now let $q$ be a prime number such that $q \nmid x$ and $q \equiv ad^{-1} \pmod b$. For each $j \in \llbracket 1,n \rrbracket$, let $t_j$ be a prime number such that $t_j \nmid x$ and $t_j \equiv p_j^{-1} \pmod b$; observe that such a prime $t_j$ exists by Euler's theorem and Dirichlet's theorem. For every $j \in \llbracket 1,n \rrbracket$, set $x_j^i = t_jp_jqd \in M_{a,b}$ for every $i \in \llbracket 1,s_j \rrbracket$. Set
	\[
		k \coloneqq \max\left\{0, \left(\sum_{j = 1}^{n} s_j\right) - \max\left\{ 1 + \left\lceil \frac{r_1 + e_1}{r_1} \right\rceil, \ldots, 1 + \left\lceil \frac{r_m + e_m}{r_m} \right\rceil \right\} \right\}.
	\] 
	It is not hard to see that $(qd)^k\prod_{j = 1}^{n} (t_jp_jqd)^{s_j}$ is a bullet of $x$ with length the upper bound in Equation~\eqref{eq: upper bound}, which concludes our argument.
\end{proof}

\begin{remark}
	As a consequence of Propositions \ref{prop: omega primality regular case} and \ref{prop: omega primality singular case}, we obtain that $\omega(M_{a,b}) = \infty$ for every ACM $M_{a,b}$.
\end{remark}


\bigskip
\section{Length Density}
\label{sec: length density}
In this section, we consider the length density of both regular ACMs and local singular ACMs. We first note the following result achieved in ~\cite{chapman2022length}, which bounds the length density using the delta set. 

\begin{lemma}[{\cite[Proposition 3.1]{chapman2022length}}]\label{Theorem: LD Delta Bound}
For a monoid $M$ and element $x \in M^{LI}$, we have 
\begin{equation*}
\frac{1}{\textup{max} \, \Delta(x)} \leq \textup{LD}(x) \leq \frac{1}{\textup{min} \, \Delta(x)}.
\end{equation*}

\end{lemma}
We now note the following upper bound on the length of a factorization.
\begin{lemma}\label{lamma_factor}
	For an irreducible $x$, let $x = a_1a_2 \dots a_n$ where $a_1a_2 \dots a_n$ is the prime factorization of $x$ on $\mathbb{N}$. We have $n \leq \phi(b)$.
\end{lemma}
\begin{proof}
	Assume there is a set $S = \{s_1, s_2, \dots s_n\}$, where $s_i = a_1a_2\dots a_i$ and $n \geq b$. Since $x \in M$,  we have $\gcd(s_i,b) = 1$. Then, by Pigeonhole Principle, there exists $i < j$ such that $s_i \equiv s_j \pmod b$. Let $t = \frac{s_j}{s_i}$. Dividing both sides of the equality by $s_i$ gives $t \equiv 1 \pmod b$, implying $t \in M$. By Theorem \ref{thm: division inclusion}, we have that $\frac{x}{t} \in M$ as well, which is a contradiction because $x$ is irreducible. 
\end{proof}
We now consider the length density of regular ACMs.
\begin{theorem}
\label{lammap3}
	Let $M_{1,b}$ be a regular ACM. Then 
	\begin{equation*}
		\textup{LD}(M_{1,b}) =  \begin{cases}
			\varnothing & \phi(b) \leq 2\\
			\frac{1}{\phi(b)-2} & \phi(b) \geq 3
		\end{cases}.
	\end{equation*}
\end{theorem}
\begin{proof}
It has been shown in \cite{banister2007arithmetic} that when $\phi(b) \leq 2$, $M$ is half-factorial and thus the length density does not exist. 
    
Now we will show $\textup{LD}(M_{1,b}) = \frac{1}{\phi(b)-2}$ for $\phi(b) \geq 3$. Note that $\textup{LD}(M_{1,b}) \leq \frac{1}{\phi(b)-2}$ by the following construction. Let $a$ be an integer with order $\phi(b) \pmod b$, and let $a^{-1}$ be its inverse $\pmod b$. Note that by Dirichlet's theorem we have a prime $a_1 \equiv a \pmod b$ and a prime $b_1 \equiv a^{-1} \pmod b$. So $a_1^{\phi(b)}b_1^{\phi(b)}$ has solely two factorizations, $(a_1b_1)^{\phi(b)}$ and $(a_1^{\phi(b)})(b_1^{\phi(b)})$. This implies $\textup{LD}(x) = \frac{1}{\phi(b)-2}$ and $\textup{LD}(M) \leq \frac{1}{\phi(b)-2}$. Assume towards a contradiction that $\textup{max} \, \Delta(x) \geq \phi(b)-1$. This implies that there exists $x \in M_{1,b}$ such that 
\begin{equation*}
		x = a_1a_2 \dots a_n = b_1b_2 \dots b_m,
\end{equation*}
where $n > m$, $n - m \geq \phi(b)-1$, and all $a_i$ and $b_i$ are irreducible. Additionally, $x$ must have no factorization of length $l$ where $n > l > m$. We will now induct on the value of $m$. 
	
When $m = 2$, we assume there exists $x \in M_{1,b}$ such that $x = a_1a_2 \dots a_n = b_1b_2$ where $n - 2 \geq \phi(b)-1$, or $n \geq \phi(b)+1$. By Lemma \ref{lamma_factor}, $b_1$ and $b_2$ have at most $\phi(b)$ prime factors. Thus, by the Pigeonhole Principle, there exists an $a_i$ that only has one prime factor. This $a_i$ must divide either $b_1$ or $b_2$, implying that one of them is not irreducible and giving us a contradiction.
	
	For $m = k$, assume there does not exist $x \in M_{1,b}$ such that $x = a_1a_2 \dots a_n = b_1b_2 \dots b_m$ where $n > m$ and $n - m \geq \phi(b)-1$, and $x$ has no factorizations of length $l$ such that $n > l > m$.
	Now, we will prove that for $m = k+1$ there also does not exist such a $x$. First, note that Lemma \ref{lamma_factor} implies that the prime factorization of $x$ contains at most $m \phi(b)$ primes. Thus, by the Pigeonhole Principle, there exists an $a_i$ that contains $m-1$ primes. Now, if we treat each prime as distinct (including those of the same value), since there are $m$ $b_j$ in total, there must exist a $b_j$ such that $a_i$ shares no primes with $b_j$. This allows us to consider the following factorization of $x$. Let $x = (b_j)(a_i)(q)$. When $q$ can be factored into more than $m-2$ irreducibles, $x = (b_j)(a_i)(q)$ has a factorization of length $l$ such that $l > m$. By our assumption on $m = k$, the maximum length of a factorization of $(a_i)(q)$ is $q+p-3$. So, $n > l$, which implies $n > l > m$, which is a contradiction.
	
	We now consider when $q$ can be factored into less than or equal to $m-2$ irreducibles. Factoring $(a_i)(q)$, we can use $b_j$ combined with primes in $q$ to produce at most $\phi(b)$ irreducibles. We can also use $q$ itself to produce at most $m-2$ irreducibles. Therefore we have $\phi(b) + m-2$ from $(a_i)(q)$. Adding on $b_j$ gives our maximum factorization length of $\phi(b)+m-1$, which is only $\phi(b)-2$ larger than $m+1$. This gives another contradiction. Since we have reached a contradiction in both cases, by applying Lemma \ref{Theorem: LD Delta Bound} we have that $\frac{1}{\phi(b)-2} \leq \textup{LD}(M_{1,b}).$ Thus, $\textup{LD}(M_{1,b}) = \frac{1}{\phi(b)-2}$. 
\end{proof}

We additionally state the length density for the the monoid $M_{b, b}$ where $b$ is an integer with more than one prime factor.
\begin{prop}
		For a monoid $M_{b,b}$ where $b = p_1^{a_1}p_2^{a_2} \dots p_n^{a_n}$ such that $p_1, p_2, \dots p_n \in \mathbb{P}$ and $a_1, a_2, \dots a_n \in \mathbb{N}$ with $n \geq 2$, we have $\textup{LD}(M_{b,b})=1$.
\end{prop}

\begin{proof}
 Note that every element in $M_{b,b}$ can be represented as $b^{k}m$ where $b$ does not divide $m$. We will prove that there exists a factorization of any length between $2$ to $k$ inclusive. Consider the following construction for a factorization of length $c$ where $c$ is an arbitrary integer between $2$ and $k$ inclusive. Let $k-(c-2) = d$. Note that 
\begin{equation*}
b^{k}m = b^{c-2}(p_1^{(d-1)a_1}p_2^{a_2} \dots p_n^{a_n}t_1)(p_1^{a_1}p_2^{(d-1)a_2}\dots p_n^{(d-1)a_n}t_2),
\end{equation*}
where $t_1$ contains all powers of $p_1$ in $m$ and $t_2$ contains all powers of $p_2,\dots , p_n$ in $m$ is a valid construction. Thus, the length density is $1$.
\end{proof}

Now we discuss $\textup{LD}(M_{a,b})$ where $\gcd(a,b) = p^\alpha$ for $p \in \mathbb{P}$ and $\alpha \in \mathbb{N}$. Let $\beta$ denote the least integer such that $p^{\beta} \in M$. Additionally, let $a' = \frac{a}{d}$, $b' = \frac{b}{d}$, and let $\delta(\alpha,\beta)$ denote the largest integer less than $\frac{\beta}{\alpha}$. Now, we let $M$ refer to the local singular ACM $M_{a,b}$. We first consider the following theorem proved in \cite{BCS08} by Baginski, Chapman, and Schaeffer. 

\begin{theorem} \label{Theorem: Delta Set of local singular monoids}
	For all local ACMs $M$, the delta set can be characterized as follows:
	\begin{equation*}
		\Delta(M) = 
		\begin{cases}
			\varnothing & \text{if}\ \alpha = \beta = 1 \\
			\{1\} & \text{if}\ \alpha = \beta > 1 \\
			[1, \frac{\beta}{\alpha}) & \text{if}\ \alpha < \beta
		\end{cases}.
	\end{equation*}
\end{theorem}

Note that the length density of a local singular ACM where $\alpha = \beta = 1$ does not exist. Also, note that by \ref{Theorem: LD Delta Bound} and \ref{Theorem: Delta Set of local singular monoids}, the length density of monoids $M_{a,b}$ when $\alpha = \beta > 1$ is $1$.  We now find the length density of the other local singular ACMs.
\begin{prop}
    For a local singular ACM $M$ where $\alpha \neq \beta$, the $\textup{LD}(M) = \frac{1}{\delta(\alpha,\beta)}.$
\end{prop}
\begin{proof}
First, we note $\frac{1}{\delta(\alpha,\beta)} \leq \textup{LD}(M)$. Consider a monoid $M_{a,b}$ with $\alpha$ and $\beta$ defined as in Theorem \ref{Theorem: Delta Set of local singular monoids}. Then, $\Delta(M) = [1, \frac{1}{\delta(\alpha,\beta)}]$. Thus, $\textup{max} \, \Delta(M) = \delta(\alpha,\beta)$ so  $\delta(\alpha,\beta) \leq \textup{LD}(M)$.

Now we show $\textup{LD}(M) \leq \frac{1}{\delta(\alpha,\beta)}$. 
Note that since $a'$ and $b'$ are relatively prime, by Dirichlet's Theorem there exists a prime $r$ such that $r \equiv a  \pmod b$. Thus, $p^{\alpha}r \in M$. Now, let $c = \delta(\alpha,\beta) \cdot(\alpha+1) - \beta $. Then we have that 
\begin{equation*}
	(p^{\alpha+c}r^k) = \frac{(p^\alpha r)^{\delta(\alpha,\beta)+2}}{(p^{\beta})}. 
\end{equation*}
Note that by Theorem \ref{thm: division inclusion} we know $(p^{\alpha+c}r^k) \in M$. So, we have
\begin{equation*}
	p^{\delta(\alpha,\beta)+2}r^{\delta(\alpha,\beta)+2} = (p^{\alpha+c}r^k)(p^{\beta}) = (p^\alpha r)^{\delta(\alpha,\beta)+2}. 
\end{equation*}
Note that there are no factorizations of $p^{\delta(\alpha,\beta)+2}r^{\delta(\alpha,\beta)+2}$ with length $l$ such that $l > \delta(\alpha,\beta)+2$ or $l < 2$. We will now prove that there are no factorizations of $p^{\delta(\alpha,\beta)+2}r^{\delta(\alpha,\beta)+2}$ of length $l$ in such that $l$ is in the interval $[3, \delta(\alpha,\beta)+1]$.

Consider an irreducible $i = p^{y_1}r^{y_2} \in M$. By Theorem \ref{Theorem: p^beta 1 mod b'}, this implies $p^{y_1}r^{y_2} \equiv 1 \pmod {b'}$. Note that $p^{\alpha y_2}r^{y_2} \equiv 1 \pmod {b'}$. Thus, $p^{\alpha y_1-y_2} \equiv 1 \pmod {b'}$. However, we know that the order of $p \pmod {b'}$ is $b$. Thus, $y_1-y_2$ is a multiple of $b$. If $y_1-y_2 = 0$, then $i$ is not an irreducible, which gives us a contradiction. However, if $y_1-y_2 \neq 0$, $\alpha y_1$ or $y_2$ must be greater than or equal to $b$, in which case the maximum length of the factorization containing $i$ is $2$. Thus, $p^{\delta(\alpha,\beta)+2}r^{\delta(\alpha,\beta)+2}$ has only factorizations of length $2$ and length $\delta(\alpha,\beta)+2$. So,
\begin{equation*}
	\textup{LD}(p^{\delta(\alpha,\beta)+2}r^{\delta(\alpha,\beta)+2}) = \frac{1}{\delta(\alpha,\beta)}.
\end{equation*} This implies $\textup{LD}(M) \leq \frac{1}{\delta(\alpha,\beta)}$.
\end{proof}

Motivated by our previous results, we conjecture that for global singular ACMs, the following holds.
\begin{conj}
For a global singular ACM M,
\begin{equation*}
    \textup{LD(}M\textup{)} = \frac{1}{\textup{max} \, \Delta(M)}.
\end{equation*}
\end{conj}

Currently, the delta set of global singular ACMs remains an open question and is likely a necessary prerequisite for determining the length density of global singular ACMs.

\bigskip
\section{Catenary Degree}
\label{sec: catenary degree}

In this section, we will determine the catenary degree of local singular ACMs. From Theorem~\ref{thm:krull}, we know that all regular ACMs are Krull. The catenary degree of Krull monoids has been bounded in~\cite{geroldinger2011catenary} and~\cite{geroldinger2015catenary}. However, computing the catenary degree of singular ACMs has remained an open problem. 

In our computation of the catenary degree, we will be making extensive use of Theorem~\ref{thm: division inclusion} to create factorization chains. Note that we can rewrite a local singular ACM $M_{a,b}$ to become $M = M_{ap^\alpha, bp^\alpha}$ for $p$ prime and $\gcd(a,b) = 1.$ Let $\beta$ be the smallest power of $p$ such that $p^{\beta} \in M.$ 

\begin{theorem}\label{thm: catenary degree}
	The catenary degree of the local singular arithmetic congruence monoid $M = M_{ap^\alpha, bp^\alpha}$ can be defined as follows:
	\begin{equation*}
		\cc(M) = 
		\begin{cases}
			2 & \text{if}\ \alpha = \beta = 1 \\
			3 & \text{if}\ \alpha = \beta >1 \\
			1 + \left\lceil \frac{\beta}{\alpha}\right\rceil & \text{if}\ \alpha < \beta
		\end{cases}.
	\end{equation*}
\end{theorem}

To begin, we will compute the catenary degree of the local singular ACM with $\alpha = \beta = 1.$ First, note the following structural information about the monoid. Since $\alpha = 1,$ we have $M = M_{ap, bp}$ for $\gcd(a,b) = 1.$ Since $\beta = 1,$ we have $p \in M.$ All elements of $M$ can be written in the form $ap + kbp$ for $k \in \mathbb{N}_0$, meaning that $p \in M$ if and only if $a = 1.$ Thus, if $\alpha = \beta =1$ then $M = M_{p, bp}$ for $b \in \mathbb{N}.$ Consider the following characterization of the set of irreducibles.

\begin{prop} 
	For $M = M_{p, bp},$ the set $$\A(M) = \{a \mid v_p(a) = 1\}.$$
\end{prop}

If we consider any element $m \in M$ with $v_p (m) >1,$ we can write $m$ as $qp^r$ in $\mathbb{N}$ for $r >1$ and $p \nmid q.$ So, $m = p \cdot qp^{r-1}$ in $\mathbb{N},$ for $r-1 >0.$ By Theorem~\ref{thm: division inclusion}, both $p$ and $qp^{r-1}$ are elements of $M,$ so any $m$ with $v_p (m)>1$ cannot be irreducible. 

Then, since all elements of $M_{p, bp}$ take the form $p + kbp$ and are thus divisible by $p,$ for all $a \in M_{p, bp}$ with $v_p (a) = 1,$ it is impossible to factor $a$ into two other elements of $M$ also with $v_p$ at least $1,$ meaning that all $a \in \{a \mid v_p(a) = 1\}$ are irreducibles. 

Now we will construct a $2$-chain $z_0, \dots, z_n$ from any factorization $z = z_0$ of an element $m = p^rq \in M$ for $p \nmid q$ to the factorization $z_n = z' = p \cdot p \cdots p \cdot pq.$ If we reorder and index the atoms within $z_i$ for $0 \leq i < n$ so that $z_i = (p\alpha_1)(p\alpha_2)\cdots(p\alpha_{r-1})(p\alpha_r)$ for $1 \leq \alpha_1 \leq \alpha_2 \leq \dots \leq \alpha_r \leq q,$ let $z_{i+1} = (p\alpha_1)(p\alpha_2)\cdots(p)(p\alpha_{r-1}\alpha_r).$ 

If both $p\alpha_{r-1}, p\alpha_{r} \in M,$ then their product $p^2\alpha_{r-1}\alpha_r \in M.$ Then, by Theorem~\ref{thm: division inclusion}, $p\alpha_{r-1}\alpha_r \in M.$ Since $v_p(p) = v_p (p\alpha_{r-1}\alpha_r) = 1,$ both factors are irreducibles. 
Note that, under this construction, the final element of this chain will be $z'.$ Since $z_{i}$ and $z_{i+1}$ differ by the relation $(p\alpha_{r-1})(p\alpha_r) = p(p\alpha_{r-1}\alpha_r),$ $\dd(z_i, z_{i+1}) = 2$ for $0 \leq i < n.$ 

So, given any two factorizations $z_1, z_2 \in \Z(m),$ we can construct a $2$-chain between $z_1$ and $z_2$ by using our construction to find a $2$-chain from $z_1$ to $z'$, then from $z'$ to $z_2,$ making $\cc(m) = 2.$ Since this construction works for all $m \in M =M_{p, bp},$ $\cc(M) = 2.$ Also, note that since each factorization of $m\in M$ must have length $v_p(m),$ $M$ is half factorial. 


Then, we will compute the catenary degree of a local singular monoid $M = M_{ap^{\alpha}, bp^{\alpha}}$ for $\gcd(a,b) = 1$ and $\alpha >1.$ Using a similar logic as in the previous section, if $p^\alpha \in M$ then $a = 1.$ So, $M$ takes the form $M_{p^\alpha, bp^\alpha}$ for $b \in \mathbb{N}$ and $p \nmid b.$ As before, we will first characterize the irreducibles of this monoid.  

\begin{prop}\label{prop:powerirr}
	For $M = M_{p^\alpha, bp^\alpha},$ the set
	$$\A(M) = \{a \mid \alpha \leq v_p (a) \leq 2\alpha-1\}.$$
\end{prop}

Suppose we have some $m \in M$ for $m = p^{n\alpha + n'}q$ for $0 \leq n' < \alpha.$ If $n > 1,$ we can factor $m$ in $\mathbb{N}$ as $m = p^\alpha \cdot p^{(n-1)\alpha + n'}q.$ Then, by Theorem~\ref{thm: division inclusion}, both $p^{\alpha}$ and $p^{(n-1)\alpha + n'}q$ are elements of $M$, so any $m = p^{n\alpha + n'}q$ with $n > 1$ is reducible. 

Also, since all elements in $M$ take the form $p^\alpha + kbp^\alpha$ and are divisible by $p^\alpha,$ for all $a \in M_{p^\alpha, bp^\alpha}$ with $\alpha \leq v_p(a) \leq 2\alpha -1,$ it is impossible to factor $a$ into two other elements of $M,$ meaning that all $a$ in the set described in~\ref{prop:powerirr} are irreducible in $M.$ 

Now we will construct a $3$-chain $z_0, \dots, z_n$ from any factorization $z = z_0$ of an element $m = p^{n\alpha + n'}q$ in $M$ to the factorization $z' = p^\alpha \cdot p^\alpha \cdots p^\alpha \cdot p^\alpha q.$ If we reorder and reindex the atoms within a factorization $z_i$ for $0 \leq i < n$ so that 
$$z_i = (p^{\alpha + e_1}\beta_1)(p^{\alpha + e_2}\beta_2)\cdots(p^{\alpha + e_{r-1}}\beta_{r-1})(p^{\alpha+e_{r}}\beta_{r})$$
for $0 \leq e_1 \leq e_2 \dots \leq e_r < \alpha$ and for $\beta_i \leq \beta_{i+1}$ if $e_i = e_{i+1},$ we define $z_{i+1}$ as follows: 
\begin{enumerate}
	\item If $0 < e_{r-1} + e_r < \alpha$ then let 
	$$z_{i+1} = (p^{\alpha + e_1}\beta_1)(p^{\alpha + e_2}\beta_2)\dots(p^{\alpha})(p^{\alpha+e_{r-1}+ e_{r}}\beta_r\beta_{r-1}).$$ 
	We have that $p^{\alpha} \in \A(M)$ and $p^{\alpha+e_{r-1}+ e_{r}}\beta_r\beta_{r-1} \in \A(M)$ by Proposition~\ref{prop:powerirr} and Theorem~\ref{thm: division inclusion}, respectively. In this case, the relation $p^{\alpha + e_{r-1}}\beta_{r-1} \cdot p^{\alpha + e_r}\beta_{r} = p^\alpha \cdot p^{\alpha + e_{r-1} + e_r}\beta_{r-1}\beta_r$ relates $z_i$ and $z_{i+1},$ meaning that $\dd(z_i,z_{i+1}) = 2.$
	
	\item If $\alpha \leq e_{r-1} + e_r < 2\alpha -1$ then 
	$$z_{i + 1} = (p^{\alpha + e_1}\beta_1)(p^{\alpha + e_2}\beta_2)\cdots(p^{\alpha})(p^{\alpha})(p^{e_{r-1}+ e_{r}- \alpha}\beta_r\beta_{r-1}).$$ 
	Again, $p^{\alpha} \in \A(M)$ by Proposition~\ref{prop:powerirr}, and $p^{e_{r-1}+ e_{r}- \alpha}\beta_r\beta_{r-1} \in \A(M)$ by Theorem~\ref{thm: division inclusion} and Proposition~\ref{prop:powerirr}. The relation $p^{\alpha + e_{r-1}}\beta_{r-1} \cdot p^{\alpha+ e_r}\beta_r = p^{\alpha}\cdot p^{\alpha}\cdot p^{e_{r-1}+ e_{r}- \alpha}\beta_r\beta_{r-1}$ relates $z_i$ and $z_{i+1},$ meaning that $\dd(z_i,z_{i+1}) = 3.$
	
	\item If $e_{r-1} + e_{r} = 0$ then both $e_{r-1}$ and $e_r$ are $0.$ Then, $z_i$ must take the form 
	$$z_i = (p^{\alpha}\beta_1)(p^{\alpha }\beta_2)\cdots(p^{\alpha}\beta_{r-1})(p^{\alpha}\beta_{r}).$$
	By Theorem~\ref{thm: division inclusion} and Proposition~\ref{prop:powerirr}, $p^{\alpha}\beta_{r-1}\beta_r \in \A(M).$ So, we use the relation $p^{\alpha}\beta_{r-1}\cdot p^{\alpha}\beta_{r} = p^{\alpha}\cdot p^{\alpha}\beta_{r-1}\beta_{r}$ to create $$z_{i+1} = (p^{\alpha}\beta_1)(p^{\alpha }\beta_2)\cdots(p^{\alpha})(p^{\alpha}\beta_{r-1}\beta_{r}),$$
	meaning that $\dd(z_{i}, z_{i+1}) = 2$ in this case. 
\end{enumerate}

Under this construction, $z_n = (p^{\alpha})\cdots (p^{\alpha})(p^{\alpha + n'}\beta_1\dots \beta_r)$ for $n' \equiv e_1 + e_2 + \dots + e_r \pmod \alpha,$ which is equal to $z',$ meaning that we have created a chain of factorizations $z_0, \dots, z_n$ from any factorization $z = z_0$ of $m \in M_{p^\alpha, bp^\alpha}$ to $z' = z_n$ such that the distance between two adjacent factorizations is no more than $3.$ Then, for any $z_1, z_2 \in \Z(m),$ we can construct a $3$-chain between them by making a $3$-chain from $z_1$ to $z'$ then $z'$ to $z_2,$ meaning that $\cc(m) = 3$ for all $m \in M.$ It follows that $\cc(M) = 3$ as well.

Now, we will compute the catenary degree of the local singular ACM with $\alpha < \beta,$ which we claim to be $1 + \left\lfloor 
\frac{\beta}{\alpha} \right\rfloor.$ However, before doing so, we will partially characterize the reducibles and irreducibles of the monoid. 

\begin{prop}
	For $M = M_{ap^\alpha, bp^\alpha}$ for $\gcd(a,b) = 1$ and $p^{\beta}$ being the smallest power of $p$ in $M,$ all numbers $m \in M$ such that $v_p (m) \geq \alpha + \beta$ are reducible. 
\end{prop}

First, consider some $m \in M$ that can be factored in $\mathbb{N}$ as $m = p^{\alpha + \beta + n}q$ for $n \in \mathbb{N}_0$ and $p \nmid q.$ Then, consider the factorization of $m$ in $\mathbb{N}$ $m = p^{\beta} \cdot p^{\alpha+n}q.$ By Theorem~\ref{thm: division inclusion}, $p^{\alpha+n}q \in M,$ so $m$ must be reducible. 

Moreover, note that all $m \in M$ must have $v_p(m) \geq \alpha,$ and that all $a \in M$ with $v_p (a) < 2\alpha$ must be irreducible, since it is impossible to write $a$ as a product of two integers both divisible by $p^\alpha.$ Then, consider the following theorem, which provide a lower bound to the catenary degree of $M$ using the delta set of $M.$ 

\begin{theorem}[{\cite[Theorem~1.6.3]{GHK06}}]
	If $M$ is a non-factorial monoid, then $2 + \sup\Delta(M) \leq \cc(M).$
\end{theorem}

Then note that by Theorem \ref{Theorem: Delta Set of local singular monoids}, $\max\Delta(M) \leq \beta/\alpha,$ so $\cc(M) > 2 + \beta/\alpha.$ This is equivalent to saying $\cc(M) \geq 1 + \lceil \beta / \alpha \rceil,$ since the catenary degree must be integral. We then provide an upper bound on $\cc(M)$.


\begin{prop}
	For $c = 1 + \left\lfloor \frac{\alpha + \beta -1}{\alpha} \right\rfloor,$ we can construct a $c$-chain between any two factorizations $z_1$ and $z_2$ of all numbers $m \in M.$ That is, $\cc (M) \leq 1 + \left\lfloor \frac{\alpha + \beta -1}{\alpha} \right\rfloor.$
\end{prop} 

Consider the element $m \in M,$ which can be factored in $\mathbb{N}$ as $p^{n\beta + k}q$ for $k < \alpha + \beta.$ We will describe a way to construct a $c$-chain for $c$ at most $1 + \left\lfloor \frac{\alpha + \beta -1}{\alpha} \right\rfloor$ between any factorization $z \in \Z(m)$ and the factorization $z' = p^{\beta} \cdot p^{\beta} \cdots p^{\beta} \cdot \prod a_i = \left(p^{\beta}\right)^n \cdot \prod a_i,$ for $a_i$ a specific product of atoms multiplying to $p^{k}q.$ 

Let the chain be $z_0, z_1, \dots, z_t$ so that $z = z_0.$ Then, if we index the atoms within the factorization $z_i = (p^{\alpha + e_1}q_1)(p^{\alpha + e_2}q_2)\cdots(p^{\alpha + e_r}q_r)$ for $0 \leq e_1 < \beta,$ $p \nmid q_i,$ and $1 \leq q_1 \leq q_2 \leq \dots \leq q_r,$ we construct $z_{i+1}$ as follows. Consider the atoms at the ``end" of the factorization of $z_i,$ starting with $p^{\alpha + e_r}q_r, p^{\alpha + e_{r-1}}q_{r-1}, \dots.$ Let $s$ be the smallest number such that $v_p\left(\prod_{i = 0}^s p^{\alpha + e_{r-s}}q_{r-s}\right) \geq \alpha + \beta.$ Note that $v_p\left(\prod_{i = 0}^s p^{\alpha + e_{r-s}}q_{r-s}\right) < 2\alpha + 2\beta -1,$ or else this would contradict $s$ being the smallest such number. Then, we have two cases.
\begin{enumerate}
	\item If $\alpha + \beta \leq v_p\left(\prod_{i = 0}^s p^{\alpha + e_{r-s}}q_{r-s}\right) < \alpha + 2\beta,$ we let 
	\[z_{i+1} = (p^{\alpha + e_1}q_1)(p^{\alpha + e_2}q_2)\cdots(p^{\alpha + e_{r-s-1}}q_{r-s-1})(p^{\beta})\prod a'_i\]
	for some set of $a'_i \in \A(M)$ with $\prod a'_i = p^{(s+1)\alpha + e_{r} + \dots +e_{r-s} -\beta}q_{r}q_{r-1}c\dots q_{r-s}.$ Then, $z_i$ and $z_{i+1}$ are related by the relation 
	\[ \prod_{i = 0}^s p^{\alpha + e_{r-s}}q_{r-s} = p^{\beta} \cdot \prod a'_i.\]
	
	There are at most $\lceil \frac{\alpha + \beta}{\alpha}\rceil$ atoms on the left-hand side of this equation, and at most $1 + \lfloor \frac{\alpha + \beta - 1}{\alpha}\rfloor$ atoms on the right of the equation. So, $\dd(z_i, z_{i+1}) \leq 1 + \lfloor \frac{\alpha + \beta - 1}{\alpha}\rfloor = \lceil \frac{\alpha + \beta}{\alpha}\rceil.$
	
	\item If $\alpha + 2\beta \leq v_p\left(\prod_{i = 0}^s p^{\alpha + e_{r-s}}q_{r-s}\right) < 2\alpha + 2\beta - 1,$ then we let 
	\[z_{i+1} = (p^{\alpha + e_1}q_1)(p^{\alpha + e_2}q_2)\cdots(p^{\alpha + e_{r-s-1}}q_{r-s-1})(p^{\beta})(p^{\beta})\prod a'_i\]
	for some set of $a'_i \in \A(M)$ with $\prod a'_i = p^{(s+1)\alpha + e_{r} + \dots e_{r-s} -2\beta}q_{r}q_{r-1}\dots q_{r-s}.$ Then, $z_i$ and $z_{i+1}$ are connected via the relation 
	\[ \prod_{i = 0}^s p^{\alpha + e_{r-s}}q_{r-s} = p^{\beta} \cdot p^{\beta} \cdot \prod a'_i.\]
	
	There are at most $\lceil \frac{\alpha + \beta}{\alpha}\rceil$ atoms on the left-hand side of this equation, and at most $2 + \lfloor \frac{2\alpha -2 1}{\alpha}\rfloor$ atoms on the right-hand side of the equation. Note that $\lceil \frac{\alpha + \beta}{\alpha}\rceil \geq 2 + \lfloor \frac{2\alpha -2 1}{\alpha}\rfloor,$ so $\dd(z_i, z_{i+1}) \leq \lceil \frac{\alpha + \beta}{\alpha}\rceil.$ 
\end{enumerate}

In both cases, $\dd(z_i, z_{i+1}) \leq \lceil \frac{\alpha + \beta}{\alpha}\rceil = 1+ \lfloor \frac{\alpha + \beta - 1}{\alpha} \rfloor.$ Then, note that under this construction, the final factorization
\[z_t = (p^\beta)(p^\beta) \cdots (p^\beta) \cdot \prod a'_i = (p^{\beta})^n \cdot \prod a'_i\]
for $\prod a'_i = p^kq.$ Then, to this chain, we add one final element $z' = \left(p^{\beta}\right)^n \cdot \prod a_i$ using the relation 
\[\prod a'_i = \prod a_i.\]
Since $\prod a'_i = \prod a_i = p^kq,$ and $k < \alpha + \beta,$ there are at most $\lfloor \frac{\alpha + \beta - 1}{\alpha}\rfloor$ atoms on both sides of the relation, meaning that $\dd(z_t, z') \leq \lfloor \frac{\alpha + \beta - 1}{\alpha}\rfloor.$ Then, the distance between any two adjacent elements in the entire chain $z_0, z_1, \dots, z_t, z'$ is bounded above by $1+ \lfloor \frac{\alpha + \beta - 1}{\alpha} \rfloor,$ so we have constructed a $c$-chain from any element $z \in \Z(m)$ to $z'$ for $c \leq 1+ \lfloor \frac{\alpha + \beta - 1}{\alpha} \rfloor.$ So, for any two $z_1, z_2 \in \Z(m)$ for $m\in M,$ we can construct a $c$-chain for $c \leq 1+ \lfloor \frac{\alpha + \beta - 1}{\alpha} \rfloor$ from $z_1$ to $z',$ then from $z'$ to $z_2$ using this method, meaning that $\cc(m) \leq 1+ \lfloor \frac{\alpha + \beta - 1}{\alpha} \rfloor$ for all $m \in M,$ or that $\cc(M) \leq 1+ \lfloor \frac{\alpha + \beta - 1}{\alpha} \rfloor.$

We have shown that 
\[1+ \lfloor \frac{\alpha + \beta - 1}{\alpha} \rfloor \geq \cc(M) \geq 1 + \lceil \frac{\beta}{\alpha}\rceil.\] 
Since $1 + \lceil \frac{\beta}{\alpha}\rceil = 1+ \lfloor \frac{\alpha + \beta - 1}{\alpha} \rfloor,$ it follows that $\cc(M) = 1 + \lceil \frac{\beta}{\alpha}\rceil$ for $M$ a local singular ACM with $\alpha < \beta.$ 
Thus, we have resolved the catenary degree for all local singular ACMs. As a corollary, consider the following observation.

\begin{cor}
For all $n \geq 2,$ we can find a local singular ACM $M$ such that $\cc(M) = n.$
\end{cor}

Consider $M=M_{p^{n-1}, (p^{n-1}-1)p},$ which has $\alpha = 1, \beta = n-1.$ Then, by Theorem~\ref{thm: catenary degree}, $\cc(M) = 1+ \left\lceil \frac{n-1}{1} \right\rceil = n.$ Besides Theorem~\ref{thm: catenary degree}, we will also propose a conjecture regarding the catenary degree of global singular ACMs. Consider the global singular ACM $M = M_{ad, df}$ such that $\gcd(a,f) =1$ and $d = \prod_{i =1}^np_i^{\alpha_i}$ for $n>1.$ 

We will define an analogous structure to $\beta$ in the local singular case. Let the set $X$ denote the set of all $x = \prod_{i=1}^n p_i^{k_i\alpha_i} \in M.$ Then, define 
\[\zeta = \min_{x\in X}\left\{\max_{1 \leq i \leq n} \left\{k_i : x = \prod_{i=1}^n p_i^{k_i\alpha_i}\right\}\right\}.\]
Let $\mu\in X$ be the element where we find $\zeta.$ Then, let $\mu'$ be the element with the second lowest 
\[\max_{1 \leq i \leq n} \left\{k_i : x = \prod_{i=1}^n p_i^{k_i\alpha_i}\right\}.\]

\begin{definition}
     Let the \textit{catenary order} $\omega_m$ of an element $m\in M$ be the least power of $m$ such that $m^{\omega_m}$ does not have a unique factorization.
\end{definition}

With these definitions, we propose the following conjecture for the catenary degree of global singular ACMs, which remains to be resolved:

\begin{conj}
For a global singular ACM $M,$ 
$$\cc(M) = \max\{\zeta + 1, \omega_{\mu},\cc\left((\mu')^{\zeta} \cdot \mu^{\omega -1}\right)\}.$$
\end{conj}

\bigskip
\section*{Acknowledgments}
 During the preparation of this paper, the last three authors were participating in PRIMES, and they would like to thank Harold Polo for the valuable feedback as well as many useful suggestions we received from him during this period. We would also like to thank the PRIMES program for giving us the opportunity to participate in pure math research and learn math topics that we would not otherwise touch upon in high school.



\end{document}